\documentclass[preprint,12pt]{elsarticle}
\usepackage{amsfonts,amssymb}
\usepackage{amsmath}
\usepackage[capitalize]{cleveref}

\newcommand{\EQ}{\begin{equation}}
\newcommand{\EN}{\end{equation}}

\newtheorem{theo}{Theorem}
\newtheorem{corollary}[theo]{Corollary}
\newtheorem{lem}[theo]{Lemma}
\newtheorem{prop}[theo]{Proposition}
\newtheorem{conjecture}[theo]{Conjecture}
\newtheorem{deff}[theo]{Definition}
\newtheorem{rem}[theo]{Remark}
\newtheorem{example}[theo]{Example}

\newcommand{\bx}{{\bf x}}

\newcommand{\bv}{{\bf v}}
\newcommand{\bu}{{\bf u}}

\newcommand{\F}{\mathbb{F}}

\newcommand{\pr}{\noindent{\bf Proof. \ }}

\newcommand{\CR}{\operatorname{CR}}
\newcommand{\UPWS}{\operatorname{UPWS}}
\newcommand{\MDS}{\operatorname{MDS}}

\journal{Finite Fields and their Applications}

\begin{document}

\begin{frontmatter}

\title{On maximum distance separable and completely regular codes}

\author[JB]{Joaquim Borges}
\author[JB]{Josep Rif\`{a}}
\author[VZ]{Victor A. Zinoviev}

\affiliation[JB]{organization={Department of Information and Communications Engineering, Universitat Aut\`{o}noma de Barcelona}}
%\affiliation[JR]{organization={Department of Information and Communications Engineering, Universitat Aut\`{o}noma de Barcelona}}
\affiliation[VZ]{organization={A.A. Kharkevich Institute for Problems of Information Transmission, Russian Academy of Sciences}}

\begin{abstract}
    We investigate when a maximum distance separable ($\MDS$) code over $\F_q$ is also completely regular ($\CR$).
For lengths $n=q+1$ and $n=q+2$ we provide a complete classification of the $\MDS$ codes that are $\CR$ or at least uniformly packed in the wide sense ($\UPWS$).
For the more restricted case $n\leq q$ with $q\leq 5$ we obtain a full classification (up to equivalence) of all nontrivial $\MDS$ codes: there are none for $q=2$; only the ternary Hamming code for $q=3$; four nontrivial families for $q=4$; and exactly six linear $\MDS$ codes for $q=5$ (three of which are $\CR$ and one admits a self-dual version).
Additionally, we close two gaps left open in a previous classification of self-dual $\CR$ codes with covering radius $\rho\leq 3$: we precisely determine over which finite fields the $\MDS$ self-dual completely regular codes with parameters $[2,1,2]_q$ and $[4,2,3]_q$ exist. \end{abstract}


\end{frontmatter}

%\linenumbers

\section{Introduction}
Let $\F_q$ denote the finite field of order $q$, where $q$ is a prime power, and let $\F_q^*=\F_q\setminus\{0\}$.
A linear \emph{$[n,k,d]_q$ code} $C$ is a $k$-dimensional subspace of $\F_q^n$ with minimum Hamming distance $d$.
For any $\bv \in \F_q^n$ its \emph{distance} to the code $C$ is $d(\bv,C)=\min\{d(\bv,\bx)\mid \bx\in C\}$. Throughout the paper we assume $0<k<n$, thereby excluding the trivial codes consisting of only the zero codeword or the whole space $\F_q^n$.

The Euclidean \emph{inner product} of two vectors $\mathbf{v},\mathbf{u}\in\F_q^n$ is
\[
\mathbf{v}\cdot\mathbf{u}=\sum_{i=1}^{n} v_i u_i\in\F_q.
\]

The (Euclidean) \emph{dual code} of a linear code $C$ is
\[
C^\perp=\{\mathbf{x}\in\F_q^n\mid \mathbf{x}\cdot\mathbf{v}=0\ \forall\mathbf{v}\in C\}.
\]

We say that $C$ is \emph{self-dual} if $C=C^\perp$. The \emph{external distance} of $C$ (see~\cite{D73}) is the number $s$ of nonzero weights in $C^\perp$.

The \emph{packing radius} of $C$ is $e=\lfloor(d-1)/2\rfloor$ and the \emph{covering radius} is
\[
\rho=\max_{\mathbf{v}\in\F_q^n}\{d(\bv,C)\}.
\]
Clearly $e\leq\rho$, with equality if and only if $C$ is perfect. It is well known that any nontrivial perfect code  (with more than two codewords)  satisfies $e\leq 3$~\cite{Tiet,ZL73}.

Linear perfect codes with $e=1$ are the \emph{Hamming codes}; they exist for length $n=(q^m-1)/(q-1)$ ($m\geq 2$), dimension $k=n-m$, and minimum distance $d=3$. Their duals are the \emph{simplex codes} with parameters $[n,m,q^{m-1}]_q$, which are equidistant~\cite[Theorem~1.8.3]{HP} (all nonzero codewords have weight $q^{m-1}$).

For any $\bv\in\F^n_q$ and any $t\in\{0,\ldots,n\}$, define $B_{\bv,t}=|\{\bx\in C\mid d(\bv,\bx)=t\}|$.

\begin{deff}[\cite{D73}]\label{CRDel}
  A code $C$ is {\em completely regular} ($\CR$) if $B_{\bv,t}$ depends only on $t$ and $d(\bv,C)$.
\end{deff}

\begin{rem}\label{cosets}
    If $C$ is linear, then Definition \ref{CRDel} is equivalent to say that  any coset $C+\bx$ depends only on its minimum weight $d(\bx,C)$.
\end{rem}
For a more detailed description and basic properties of $\CR$ codes over finite fields, see \cite{BRZ19,previ}, or \cite[Chap. 2]{KKM}, for example.

A more general class of codes is the class of uniformly packed codes in the wide sense \cite{BZZ74}.

\begin{deff}[\cite{BZZ74}]
  A code $C\subseteq\F_q^n$ with covering radius $\rho$ is {\em uniformly packed in the wide sense} ($\UPWS$) if there exist rational numbers $\beta_0,\ldots,\beta_\rho$ such that
\begin{equation}\label{SumUPWS}
\sum_{i=0}^{\rho} \beta_i B_{\bx,i} = 1,
\end{equation}
for any $\bx\in\F^n_q$. The numbers $\beta_0,\ldots,\beta_\rho$ are called the {\em packing coefficients}.
\end{deff}

If $C$ is a $[n,k,d]_q$ code, the Singleton bound \cite{Sing} establishes that
\begin{equation}\label{singleton}
  d\leq n-k+1.
\end{equation}
If $C$ attains equality, then $C$ is a {\em maximum distance separable} ($\MDS$) code.

The Griesmer bound \cite{Grie,Grie2} is
\begin{equation}\label{griesmer}
  n\geq \sum_{i=0}^{k-1} \left\lceil \frac{d}{q^i} \right\rceil\;\;\; (k\geq 1).
\end{equation}
If $C$ attains equality, then $C$ is a {\em Griesmer} code.

In this paper we study the possible parameters of $\MDS$ codes (which are also Griesmer codes) over $\F_q$. For length $n=q+1$ and $n=q+2$ we give a complete classification of the $\MDS$ codes that are completely regular or at least uniformly packed in the wide sense. For $q\leq 5$, we obtain a full classification (up to equivalence) of all nontrivial $\MDS$ codes. For $q=2$, as is well known, there are only trivial cases. For $q=3$, there is only one nontrivial case, which is a self-dual Hamming code. For $q=4$, the nontrivial cases are a Hamming code, its dual code, a self-dual code, and the so-called hexacode. The case $q=5$ being the most interesting. We determine that there are $6$ nontrivial such codes in $\F_5$, up to equivalence. Three of them are $\CR$ codes. In addition, we point out in which cases such codes are self-dual. Moreover, we specify exactly in which finite fields the $\MDS$ codes with parameters $[2,1,2]_q$ and $[4,2,3]_q$ are also self-dual. Such codes are completely regular as stated in \cite{previ}.

The main conjecture on $\MDS$ codes is the following
\begin{conjecture}\label{MDSConj}
  If $C$ is a nontrivial $\MDS$ $[n,k,d]_q$ code ($2\leq k\leq n-2$), then $n\leq q+1$ unless $q$ is even and $k=3$ or $k=q-1$, in which case $n\leq q+2$.
\end{conjecture}

Such conjecture was first considered by Segre \cite{Seg} and is proven for a large number of cases. For our purposes, the following case, proven in \cite{Ball} is useful.

\begin{prop}\label{conj}
If $C$ is an $\MDS$ $[n,k,d]_p$ code with $p$ prime and $2\leq k\leq n-2$, then $n\leq p+1$.
\end{prop}

The paper is organized as follows.
In Section~\ref{Basic} we collect the basic definitions and known results that will be used throughout the paper, including some simple lemmas on equations over finite fields that allow us to precise, for the first time, the exact fields where self-dual $[4,2,3]_q$ $\MDS$ codes exist.

Section~\ref{trivial} deals with the trivial $\MDS$ codes ($k=1$ or $k=n-1$).
We recall their covering radius and external distance and clarify when they are completely regular or uniformly packed.

In Section~\ref{nontrivial} we study nontrivial $\MDS$ codes ($2\leq k\leq n-2$).
We first analyse the extremal cases allowed by the $\MDS$ conjecture: length $n=q+2$ (Subsection~\ref{sec:41}) and length $n=q+1$ (Subsection~\ref{sec:42}).
For $n=q+1$ we give a sharp characterisation of the parameters for which the code is uniformly packed in the wide sense (and hence, in some cases, completely regular), including the infinite family of completely regular doubly-extended Reed--Solomon codes of dimension $q-2$ over even-order fields.
Subsection~\ref{sec:43} is devoted to the short nontrivial $\MDS$ codes ($n\leq q$).
We prove that all such codes are uniformly packed and, for the smallest fields $q\leq 5$, we obtain a complete classification (Theorem~\ref{principal} for $q=5$ being the most interesting case).

Section~\ref{sec:5} summarises all families and sporadic examples obtained.
%of $\MDS$ completely regular codes
%obtained and closes the two pending self-dual cases.
% completely regular $\MDS$ codes left open in the earlier classification of self-dual completely regular codes with covering radius at most $3$.

\section{Preliminary results}\label{Basic}

\subsection{Quadratic equations over finite fields}

The next result can be derived from Remark 5.13 in \cite{LN}. However, we give here our own proof, which is very simple.

\begin{prop}\label{quadrat}
Let $q$ be an odd prime power. Then $-1$ is a square in $\F_q$ if and only if $q\equiv 1 \pmod{4}$.
\end{prop}

\pr
Let $\xi$ be a primitive element in $\F_q$. Since $\xi^{q-1}=1$, we have $\xi^\frac{q-1}{2} =-1$. But $-1$ is a square if and only if there is some element $\xi^i\in\F_q$ such that $\xi^{2i}=-1$. Therefore, $\frac{q-1}{2}$ is even, implying $q-1\equiv 0 \pmod{4}$, or the same, $q\equiv 1 \pmod{4}$.
\qed

\begin{rem}
If $q=p^r$, where $p$ is an odd prime and $r>0$, then the condition $q\equiv 1 \pmod{4}$ is clearly equivalent to:
\begin{itemize}
  \item[(i)] $p\equiv 1 \pmod{4}$, or
  \item[(ii)] $p\equiv 3 \pmod{4}$ and $r$ is even.
\end{itemize}
\end{rem}

Recall that the quadratic character $\chi$ of $\F_q$ is defined as $\chi(a)=1$ if $a$ is a square, and $\chi(a)=-1$ if $a$ is not ($a\in\F_q^*$).
The next result is Lemma 6.24 in \cite{LN}.

\begin{prop}\label{equacions}
Let $q$ be an odd prime power and $t, a_1,a_2\in\F_q^*$.

The equation
\begin{equation}\label{EqBase}
  a_1x_1^2 + a_2x_2^2 =t
\end{equation}
has exactly $q-\chi(-a_1a_2)$ solutions, where $\chi$ is the quadratic character of $\F_q$.
\end{prop}

In \cite{previ} it is stated that any $[4,2,3]_q$ code is $\CR$ and taking a generator matrix of the form
\begin{equation}\label{matriu}
G~=~\left(
\begin{array}{cccc}
~1~&~0~&~\alpha~&~\beta\\
~0~&~\xi~&~\beta~&~-\alpha\\
\end{array}
\right),
\end{equation}
where $\alpha, \beta \in \F^*_q$, such code is also self-dual if $1+\alpha^2+~\beta^2 = 0$ and $\xi^2=1$. But in \cite{previ} it is not specified in which finite fields these conditions can be satisfied. By using Proposition \ref{equacions}, we can precise such finite fields.

\begin{corollary}\label{423}
  A self-dual $[4,2,3]_q$ code exists if and only if $q\notin\{2,5\}$.
\end{corollary}

\pr
Consider, in $\F_q$, the equation
\begin{equation}\label{EqBona}
  \alpha^2+\beta^2=-1.
\end{equation}
If $q=2$, then it is clear that there is no solution such that $\alpha,\beta\neq 0$. For $q=2^r$ ($r>1$), we always have solutions with $\alpha,\beta\neq 0$, since $\alpha^2+\beta^2=(\alpha+\beta)^2$. Hence, taking any $\alpha\notin\{0,1\}$ we have the solution $(\alpha, \beta=1+\alpha)$, where indeed $\alpha,\beta\neq 0$.

Assume now that $q$ is odd. By Proposition \ref{equacions}, Eq. (\ref{EqBona}) has $q-\chi(-1)$ solutions. However, the set of solutions can include cases where $\alpha$ or $\beta$ is $0$.

If $-1$ is not a square then $\alpha$ and $\beta$ must be nonzero elements and, by Proposition \ref{equacions} the number of solutions of Eq. (\ref{EqBona}) is $q+1$. Hence, in this case we always have solutions (with $\alpha,\beta\neq 0$).

If $-1$ is a square, then let $\gamma\in\F_q^*$ be such that $\gamma^2=-1$. Thus, we have four solutions where $\alpha$ or $\beta$ is $0$:
$$
(\alpha,\beta)\in\{(\gamma,0), (-\gamma,0), (0,\gamma), (0,-\gamma)\}.
$$
On the other hand, in this case, the total number of solutions is $q-1$, according to Proposition \ref{equacions}. Therefore, Eq. (\ref{EqBona}) always has solutions with $\alpha,\beta\neq 0$ if $q-1 > 4$, that is, $q>5$.

Note that $-1$ is not a square in $\F_3$ and it is a square in $\F_5$. We conclude that Eq. (\ref{EqBona}) always has solutions with $\alpha,\beta\neq 0$ except in $\F_2$ and $\F_5$.
\qed

\subsection{Properties of completely regular and Hamming codes}

For completely regular codes, we will use the following properties.
\begin{lem}\label{CRCodes}
Let $C$ be a linear code with minimum distance $d$, covering radius $\rho$, packing radius $e$ and external distance $s$.
\begin{itemize}
  \item[(i)] $\rho\leq s$ \cite{D73}.
  \item[(ii)] If $d\geq 2s-1$, then $C$ is $\CR$ \cite{D73}.
  \item[(iii)] $\rho = s$ if and only if $C$ is $\UPWS$ \cite{BZZ77}.
  \item[(iv)] If $C$ is $\CR$, then $\rho=s$ \cite{Sole}.
  \item[(v)] If $C$ is $\UPWS$ and $\rho=e+1$, then $C$ is $\CR$ \cite{GvT,SZZ71}.
  \item[(vi)] Let $C_1$ and $C_2$ be two cosets of $C$ with the same minimum weight $w$. If $w\leq d-s$ or $w=s$, then $C_1$ and $C_2$ have identical weight distribution \cite{D73} (see also \cite[Thm. 4.3]{Sole}).
\end{itemize}
\end{lem}

Note that, by (iii) and (iv), a $\CR$ code is always a $\UPWS$ code.

The only self-dual perfect code is the ternary Hamming $[4,2,3]_3$ code (see, for example, \cite[Lemma 13]{previ}). Since a Hamming code has parameters $[n=(q^m-1)/(q-1),n-m,3]_q$, it is an $\MDS$ code if and only if $m=2$. In such case, the parameters are $[q+1,q-1,3]_q$ and the dual simplex code has parameters $[q+1,2,q]_q$. Since a simplex code has only one nonzero weight, Hamming codes have external distance $s=1$. Hence a Hamming code is $\CR$ by Lemma \ref{CRCodes} (ii).

\subsection{$\MDS$ and Griesmer codes}

The following equivalence is easy to see, but not usually mentioned.

\begin{prop}\label{MDSGriesmer}
  Let $C$ be an $[n,k,d]_q$ code with $k\geq 2$. Then $C$ is an $\MDS$ code if and only if $C$ is a Griesmer code and $d\leq q$.
\end{prop}

\pr
If $C$ is an $\MDS$ code, then $n=d+k-1$. By the Griesmer bound $(\ref{griesmer})$ we have
$$
d+k-1\geq d+\sum_{i=1}^{k-1}\left\lceil \frac{d}{q^i} \right\rceil\;\;\;\Longrightarrow\;\;\;\left\lceil \frac{d}{q^i} \right\rceil =1\;\;\forall i=1\ldots k-1.
$$
Thus $d\leq q$ and the bound is attained.

Conversely, if $d\leq q$ and $C$ is a Griesmer code, then
$$
n=\sum_{i=1}^{k-1}\left\lceil \frac{d}{q^i} \right\rceil=d+k-1.
$$
Hence, $C$ is an $\MDS$ code.
\qed

%The next result can be found in \cite{Ward}.
%begin{lem}\label{divisible}
% Let $C$ be a $[n,k,d]_p$ Griesmer code, where $p$ is a prime number. If $p^e$ divides $d$ ($e\geq 1$), then $p^e$ divides the weight of any codeword in $C$.
%\end{lem}

%Codes with parameters $[n,1,n]_q$ and $[n,n-1,2]_q$ are $\MDS$ for any $n$ and any $q$. An $[n,1,n]_q$ code is clearly CR for $q=2$ but not in general. By Lemma \ref{CRCodes} (ii), any $[n,n-1,2]_q$ is CR. For $n=2$, we have codes with parameters $[2,1,2]_q$. Such CR codes are self-dual if the generator matrix is of the form $(1\;\;\alpha)$, where $\alpha^2=-1$. In fact, as can be seen in \cite[Thm. 36]{previ} these are the only self-dual CR codes with minimum distance $d\leq 2$. Note that by Proposition \ref{quadrat} such codes exist if and only if $q\equiv 1 \pmod{4}$. Such condition is not specified in \cite{previ}.
%
%The $\MDS$ codes in the preceding paragraph are usually called {\em trivial}. Therefore, nontrivial $\MDS$ codes are those with $2\leq k\leq n-2$.

All the following results on $\MDS$ codes can be found in \cite[Sect. 7.4]{HP}.

\begin{prop}\label{MDSCodes}
Let $C$ be an $\MDS$ $[n,k,d]_q$.
\begin{itemize}
  \item[(i)] $C^\perp$ is $\MDS$.
  \item[(ii)] The number of minimum weight codewords in $C$ (i.e. codewords of weight $d$) is
  $$
  A_d=\binom{n}{d} (q-1).
  $$
  \item[(iii)] If $2\leq k\leq n-2$, then $k\leq \min\{n-2,q-1\}$ and $n\leq 2(q-1)$.
\end{itemize}
\end{prop}

In \cite{Sole2}, the number of weights of $\MDS$ codes is determined.

\begin{prop}[\cite{Sole2}]\label{numpesos}
Let $C$ be an $\MDS$ $[n,k,d]_q$ code and let $s'$ be the number of nonzero weights of $C$.
\begin{itemize}
  \item[(i)] If $C$ is an $[n,n-1,2]_2$ code with $n>2$, then $s' = \lfloor n/2 \rfloor \not= k$.
  \item[(ii)] If $C$ is a $[q+1,2,q]_q$ code, then $s'=1$.
  \item[(iii)] If $C$ is a $[2^m+2,3,2^m]_{2^m}$ code, then $s'=2$.
\end{itemize}
In all other cases $s'=k$, assuming that Conjecture \ref{MDSConj} is true.
\end{prop}

\section{Trivial $\MDS$ codes}\label{trivial}
Trivial $\MDS$ $[n,k,d]_q$ codes are those with $k=1$ or $k\geq n-1$.

For $k=1$, we have the $\MDS$ repetition codes $C_R$ with parameters $[n,1,n]_q$. Such codes and their duals, the $[n,n-1,2]_q$ codes, are $\CR$ (by Lemma \ref{CRCodes} (ii), since the codes $C_R$ have just one nonzero weight). For $q=2$, $C_R$ codes have covering radius $\rho=\lfloor n/2 \rfloor$ (see \cite[Thm. 1]{BD}) and identical external distance $s=\lfloor n/2 \rfloor$ by Proposition \ref{numpesos}. Hence, binary $C_R$ codes are $\UPWS$ by Lemma \ref{CRCodes} (iii). Moreover, since the packing radius is $e=\lfloor (n-1)/2 \rfloor$, we have that $e=\rho$ if $n$ is odd, and $e+1=\rho$ if $n$ is even. Therefore, $C_R$ is a perfect code if $n$ is odd, and a quasi-perfect $\UPWS$ code if $n$ is even. In any case, $C_R$ is a $\CR$ code by Lemma \ref{CRCodes}.

For $q>2$ and $n>2$, we have $\rho=n-\lceil n/q \rceil$ (see \cite[Thm. 1]{BD}) and the external distance is $s=n-1$ (Proposition \ref{numpesos}). Hence, $C_R$ is $\UPWS$ if and only if $q\geq n$.

\begin{example}
The ternary code $C=\{0000,1111,2222\}$ is not $\CR$: the vectors $\bv=0011$ and $\bu=1200$ are at distance 2 from $C$, but their distances to the codewords are $2,2,4$ and $2,3,3$, respectively.
\end{example}

Note that if a $C_R$ code is self-dual then $n=2$. In such case, the codes $[2,1,2]_q$ are $\CR$: any vector $\bx$ is at distance 0 or 1 from $C$. It is esay to see that, if $\bx$ is at distance 0, then $B_{\bx,0}=1$, $B_{\bx,1}=0$, and $B_{\bx,2}=q-1$. If $\bx$ is at distance 1, then $B_{\bx,0}=0$, $B_{\bx,1}=2$, and $B_{\bx,2}=q-2$. Indeed, take a generator matrix $G=(1\;\;\alpha)$ for a $[2,1,2]_q$ $C_R$ code. Then $C_R$ is self-dual if and only if $\alpha^2=-1$. In fact, these codes or direct sums of them are the only self-dual $\CR$ codes with $d\leq 2$ (see \cite[Thm. 36(i)]{previ}). However, in \cite{previ} the finite fields where such codes exist are not specified. As a direct consequence of Proposition \ref{quadrat} we can establish the following result.

\begin{theo}\label{main}
  Let $C$ be a self-dual $\CR$ $[n,k,d]_q$ code. If $d\leq 2$, then $C$ is necessarily an $\MDS$ $[2,1,2]_q$ code, or a direct sum of $j$ codes with these parameters (and $C$ is not $\MDS$ for $j>1$). Such codes exist if and only if $q$ is even or $q \equiv 1 \pmod{4}$.
\end{theo}

\section{Nontrivial $\MDS$ codes}\label{nontrivial}

%The $\MDS$ codes of the previous section are usually referred to as trivial $\MDS$ codes. Thus,
Nontrivial $\MDS$ $[n,k,d]_q$ codes are those with $2\leq k\leq n-2$. Assuming that the {\em $\MDS$ conjecture} (Conjecture \ref{MDSConj}) is true, we have that $n\leq q+1$ except in some exceptional cases where $n=q+2$. We start with these extremal cases.

\subsection{$\MDS$ codes with length $n=q+2$}\label{sec:41}

According to Conjecture \ref{MDSConj}, a nontrivial $\MDS$ code of length $n=q+2$ is a $[2^m+2,3,2^m]_{2^m}$ code or a $[2^m+2,2^m-1,4]_{2^m}$ code, where $m>1$.

Concerning complete regularity, we introduce the following result for such codes.

\begin{theo}\label{casqmesdos}
Let $C$ be an $\MDS$ $[q+2,k,d]_q$ code, where $q=2^m$ and $m>1$.
  \begin{itemize}
    \item[(i)] If $C$ is an $\MDS$ $[2^m+2,3,2^m]_{2^m}$ code with odd $m$, then $C$ is not $\UPWS$ (and so, not $\CR$). %For $m=2$, $C$ is $\CR$.
    \item[(ii)] If $C$ is an $\MDS$ $[2^m+2,2^m-1,4]_{2^m}$ code, then $C$ is $\CR$.
  \end{itemize}
\end{theo}

\pr
(i) By Proposition \ref{numpesos}, $C$ has external distance $s=2^m-1$. For $m$ odd, the covering radius of $C$ is $\rho=2^m-2$ by \cite[Thm. 4]{GK}. Hence, by Lemma \ref{CRCodes} (iii), $C$ is not $\UPWS$.

%If $m=2$, then we have the well-known two-weight $[6,3,4]_4$ hexacode (see, e.g., \cite[p. 383]{HP}). Such code is $\CR$ by Lemma \ref{CRCodes} (ii) since $s=2$. In fact, it is a particular case of (ii).

(ii) In this case, by Proposition \ref{numpesos}, $C$ has external distance $s=2$ and thus $C$ is $\CR$ by Lemma \ref{CRCodes} (ii).
\qed

\subsection{$\MDS$ codes with length $n=q+1$}\label{sec:42}

For $n=q+1$ we have $\MDS$ codes with parameters $[q+1,k,q-k+2]_q$. We start with the particular case $k=2$.

\begin{prop}\label{k2}
Let $C$ be an $\MDS$ $[q+1,2,q]_q$ code
\begin{itemize}
  \item[(i)] If $q$ is odd, then $C$ is $\UPWS$ and $\CR$ if and only if $q=3$.
  \item[(ii)] If $q$ is even, then $C$ is $\UPWS$. For $q=2^2$, $C$ is $\CR$. For $q=2^3$, $C$ is not $\CR$.
\end{itemize}
\end{prop}

\pr
The $[q+1,2,q]_q$ code $C$ is a simplex code since its dual is a $[q+1,q-1,3]_q$ Hamming code.

(i) For odd $q$, the covering radius of $C$ is $\rho=q-2=d-2$ (see \cite[Thm. 10]{BD}) and the external distance is $s=q-1=d-1$ except for $q=3$ (see Proposition \ref{numpesos} and apply it to the dual code). Hence, for $q\neq 3$, $C$ is not $\UPWS$ by Lemma \ref{CRCodes} (iii). For $q=3$, $C$ is the self-dual Hamming $[4,2,3]_3$ code which is perfect and so $\CR$.

(ii) For even $q$, we have $s=d-1$ (again applying Proposition \ref{numpesos}) and $\rho=d-1$ \cite[Rem. 2]{GK}. Hence $C$ is an $\UPWS$ code.
If $q=2^2$, the $[5,2,4]_4$ code is $\CR$ (it is a particular case of Theorem \ref{qumesu} (iii)). For $q=2^3$, we have computationally checked that a $[9,2,8]_8$ code is not $\CR$.
\qed

Now, we have the following strong result for any $k\geq 2$.

\begin{theo}\label{qumesu}
  If $C$ is a nontrivial $\MDS$ $[q+1,k,q-k+2]_q$ code, then $C$ is not $\UPWS$ (and so not $\CR$) with the following exceptions:
  \begin{itemize}
    \item[(i)] $k=q-1$, in which case $C$ is a Hamming perfect and $\CR$ code.
    \item[(ii)] $k=2$ with $q$ even or $q=3$, in which case $C$ is $\UPWS$. For $q=3$, $C$ is also $\CR$.
    \item[(iii)] $k=q-2$ with $q$ even, in which case $C$ is $\UPWS$. Moreover, for any $m>1$, there exists a $\CR$ $[q+1,q-2,4]_q$ code, where $q=2^m$.
  \end{itemize}
\end{theo}

\pr
(i) Trivial.

(ii) Directly from Proposition \ref{k2}.

(iii) In this case, $C$ is a $[2^m+1,2^m-2,4]_{2^m}$ code. By Proposition \ref{numpesos}, the external distance of $C$ is $s=3=d-1$ and, by \cite[Rem. 2]{GK}, the covering radius of $C$ is $\rho=d-1=3$. Hence, $C$ is $\UPWS$.

Generalized Reed--Solomon (GRS) codes and their extensions are $\MDS$, in general. No $\MDS$ codes with parameters other than those arising from GRS codes or their extensions are presently known~\cite{HP}. Hence, almost all $[q+1, q-2, 4]_q$ maximum distance separable codes over $\F_q$ are equivalent to doubly-extended Reed--Solomon codes.

By Lemma \ref{CRCodes} (vi), the cosets of minimum weight $1$ have identical weight distribution. Again, by Lemma \ref{CRCodes} (vi), the cosets of weight $3$ have also the same weight distribution.

The automorphism group of these doubly-extended Reed--Solomon codes is $G=P\Gamma\mathrm{L}_2(q)$ \cite{Dur87}. Consider the action of $G$ over the cosets of $C$: for any $\sigma\in G$ and any $\bx\in\F_q^n$, $\sigma(C+\bx)=C+\sigma(\bx)$.

Since $G$ acts $2$-transitively on the $q+1$ points of the projective line over $\F_q$ \cite{DM96}, we conclude that the cosets of weight $2$ are in the same $G$-orbit and hence they have identical weight distribution.

Therefore, $C$ is $\CR$ (see Remark \ref{cosets}).

Assume now that $C$ is not one of the cases (i)--(iii). Then, $s=q-k+1=d-1$ by Proposition \ref{numpesos}. The covering radius of $C$ is $\rho=d-2$ (see \cite[Thm. 2]{GK}). Therefore, $C$ is not $\UPWS$ by Lemma \ref{CRCodes} (iii).
\qed

\begin{rem}
	The value $q = 2^2$ also shows that in Proposition~12~(ii) there exists a case where the code is completely regular.
\end{rem}

% \begin{rem}\cite[\S 11.2]{BCN89}
% 	Kasami codes are binary codes of length $n = 2^m - 1$ that can be expressed as cyclic codes $C_{1,s}$, generated by $m_1(x)m_s(x)$, where $m_i(x)$ is the minimal polynomial of $\alpha^i$ and $\alpha$ is a primitive element of $\F_{2^m}$.
	
% 	When $m$ is even, $m = 2m'$ and $s = q + 1$ with $q = 2^{m'}$, the binary Kasami code $C_{1,s}$ has parameters $[q^2 - 1, q^2 - 1 - 3m', 3]_2$; it is a \emph{completely regular} code .

%     %with intersection array
% 	%\[
% 	%\IA = [q^2-1, q^2-q, 1; \, 1, q, q^2-1].
% 	%\]
% \end{rem}

\subsection{$\MDS$ codes with length $n\leq q\leq 5$}\label{sec:43}

We begin establishing the following result.

\begin{theo}\label{enequ}
If $C$ is an $\MDS$ $[n,k,d]_q$ code with $n\leq q$, then $C$ is $\UPWS$.
\end{theo}

\pr
By Proposition \ref{numpesos}, $s=n-k=d-1$. The covering radius is $\rho=d-1$ \cite[Thm. 1]{GK}. Since $s=\rho$, $C$ is $\UPWS$ by Lemma \ref{CRCodes} (iii).
\qed

\begin{rem}
 The result $\rho=d-1$ in \cite[Thm. 1]{GK} is true, provided that $C$ can be embedded in an $[n,k+1,n-k]_q$ code, which is true for any known nontrivial $\MDS$ code of length $n\leq q$.
\end{rem}

In general, it seems hard to determine in which cases an $\MDS$ code $C$ is also $\CR$. However, for small fields we may consider some particular cases.

It is well known that all binary $\MDS$ codes are trivial. This is a direct consequence of the condition $d\leq q$ for nontrivial $\MDS$ codes (see Proposition \ref{MDSGriesmer}).

For $q=3$, we have $n\leq 4$, by Proposition \ref{MDSCodes} (iii). Hence, we have the only nontrivial case of $[4,2,3]_3$ codes which are the self-dual $\CR$ Hamming ternary codes.

If $q=4$ then, by Proposition \ref{MDSCodes} (iii), we have $n\leq 6$. So, we have the possible parameters $[6,4,3]_4$, $[6,3,4]_4$, $[6,2,5]_4$, $[5,3,3]_4$, $[5,2,4]_4$ and $[4,2,3]_4$. The case of a $[6,3,4]_4$ code corresponds to the well-known two-weight {\em hexacode} (see, e.g., \cite[p. 383]{HP}). Such code is $\CR$ by Lemma \ref{CRCodes} (ii) since $s=2$. In fact, it is a particular case of Theorem \ref{casqmesdos} (ii). The other two cases with $n=6$ can be discarded after the following result.

\begin{lem}\label{No64}
  There do not exist codes with parameters $[6,4,3]_4$ or $[6,2,5]_4$.
\end{lem}

\pr
It is enough to prove the non existence of a $[6,2,5]_4$ code since it implies the non existence of its dual, a $[6,4,3]_4$ code. Assume that $C$ is a $[6,2,5]_4$ code. Since it would be $\MDS$, by Proposition \ref{MDSCodes} (ii), the number of codewords of weight 5 should be
$$
A_5=3\binom{6}{5}=18.
$$
However, $|C|=4^2=16$ leading to a contradiction.
\qed

%\begin{rem}
%  In fact, Conjecture \ref{MDSConj} is proven for values of %$q\leq 47$ with a few exceptions. Thus, we could directly %discard the cases in Lemma \ref{No64}.
%\end{rem}

For the remaining cases, note that a $[5,3,3]_4$ code is a Hamming code and a $[5,2,4]_4$ code is the corresponding simplex dual. Finally, a $[4,2,3]_4$ code exists and it is $\CR$. In fact, any $[4,2,3]_q$ code is $\CR$ by Lemma \ref{CRCodes} since $s\leq 2$ (see also Proposition 23 in \cite{previ}). By Corollary \ref{423} a $[4,2,3]_4$ code is self-dual taking a generator matrix
$$
G=\left(
    \begin{array}{cccc}
      1 & 0 & \alpha & \alpha^2 \\
      0 & 1 & \alpha^2 & \alpha \\
    \end{array}
  \right),
$$
where $\alpha$ is a primitive element in $\F_4$.

%From the preceding paragraphs we deduce that
The next smallest interesting finite field where to study $\MDS$ codes is $\F_5$.

Now, we determine the possible parameters for nontrivial $\MDS$ codes over $\F_5$. Since Conjecture \ref{MDSConj} is proven for prime fields, and applying Proposition \ref{MDSCodes}, we only have to consider the cases:
$$
(n,k)\in\{(6,4),(6,3),(6,2),(5,3),(5,2),(4,2)\}.
$$

%\begin{lem}\label{No7No8}
%  There do not exist codes with parameters $[8,4,5]_5$, $[7,4,4]_5$, $[7,3,5]_5$.
%\end{lem}
%
%\pr
%For the case of a $[8,4,5]_5$ code $C$, by Lemma \ref{divisible}, $C$ contains only codewords of weight $5$ and the zero codeword. Hence, by (ii) in Proposition \ref{MDSCodes} the size of $C$ is
%$$
%|C|=A_5+1=4\binom{8}{5}+1=224.
%$$
%A contradiction, since the dimension is $4$ and so $|C|=5^4=625$.
%
%If a $[7,4,4]_5$ code exists then, by (i) in Proposition \ref{MDSCodes}, its $[7,3,5]_5$ dual code must also exist, as well. However, by Proposition \ref{conj}, we have that $n\leq 6$ if $2\leq k\leq 3$. Hence, the codes $[7,3,5]_5$ and $[7,4,4]_5$ do not exist.
%\qed

Finally, we establish the main result for $\MDS$ codes over $\F_5$.

\begin{theo}\label{principal}
  If $C$ is a nontrivial $\MDS$ $[n,k,d]_5$ code, then $C$ is one of the following codes, up to equivalence:
  \begin{itemize}
    \item[(i)] A $\CR$ $[6,4,3]_5$ code or its dual $[6,2,5]_5$ which is not $\UPWS$.
    \item[(ii)] A $[6,3,4]_5$ code, which is not $\UPWS$. A self-dual version of such code exists.
    \item[(iii)] A $\CR$ $[5,3,3]_5$ code or its dual $[5,2,4]_5$ which is $\UPWS$.
    \item[(iv)] A $[4,2,3]_5$ code which is $\CR$ and not self-dual.
  \end{itemize}
\end{theo}

\pr
Let us see that all these codes exist.

(i) A $[6,4,3]_5$ code is a Hamming code, which is a perfect code and hence $\CR$. The $[6,2,5]_5$ is the corresponding simplex dual code, which is not $\UPWS$ by Proposition \ref{k2}.

(ii) A $[6,3,4]_5$ code cannot be $\UPWS$ by Theorem \ref{qumesu}. A self-dual code with these parameters exists by \cite[Thm.1(iv) or (v)]{FF}.

(iii) For the $[5,2,4]_5$ code we can consider a generator matrix
$$
G=\left(
    \begin{array}{ccccc}
      1 & 0 & 1 & 1 & 1 \\
      0 & 1 & 2 & 3 & 4 \\
    \end{array}
  \right).
$$
Note that the code generated by $G$ is indeed a $[5,2,4]_5$ code, which is $\UPWS$ by Theorem \ref{enequ}. The corresponding dual code is a $[5,3,3]_5$ code which is $\CR$ by Lemma \ref{CRCodes}

(iv) Consider the generator matrix
$$
G=\left(
  \begin{array}{cccc}
    1 & 0 & 1 & 1 \\
    0 & 1 & 2 & 3 \\
  \end{array}
\right).
$$
It is readily verified that the generated code is a $[4,2,3]_5$ code $C$. The dual code $C^\perp$ has the same parameters but $C\neq C^\perp$, i.e., $C$ is not self-dual. This is clear due to Corollary \ref{423}. But the code is $\CR$ by Lemma \ref{CRCodes}. Indeed, $C$ and $C^\perp$ have external distance $s=2$ by Proposition \ref{numpesos}. In fact, $C$ and $C^\perp$ have the same weight distribution, which is:
$$
A_3=4\binom{4}{3}=16\;\;\mbox{ by Proposition \ref{MDSCodes} (ii), and } A_4=5^2-A_3-1=8.
$$
\qed

\section{Summary}\label{sec:5}
In summary we have found the following infinite families of $\MDS$ codes which are $\CR$ codes:
\begin{itemize}
  \item[(i)] The trivial $[n,n-1,2]_q$ codes for any $n$ and $q$, and the $[n,1,n]_2$ codes.
  \item[(ii)] Any $[2^m+2,2^m-1,4]_{2^m}$ code for $m>0$.
  \item[(iii)] The Hamming $[q+1,q-1,3]_q$ codes for any $q$.
  \item[(iv)] The $[2^m+1,2^m-2,4]_{2^m}$ codes equivalent to doubly-extended Reed-Solomon codes, for any $m>1$.
\end{itemize}

In addition, for $n\leq q$, we have the sporadic codes (not included above) with parameters
\begin{itemize}
    \item $[4,2,3]_4$. The punctured code of a $[5,2,4]_4$ CR code (item (iv) above, for $m=2$),
    \item $[5,3,3]_5$. Theorem \ref{principal} (iii).
    \item $[4,2,3]_5$. Theorem \ref{principal} (iv).
\end{itemize}

\section*{Acknowledgements}
This work has been partially supported by the Spanish Ministerio de Ciencia e Innovación  under Grants PID2022-137924NB-I00
(AEI/FEDER UE), RED2022-134306-T, and also by the Catalan AGAUR Grant 2021SGR-00643.

The  research of the third author of the
paper was carried out at the Institute for Information
Transmission Problems of the Russian Academy of Sciences
within
the program of fundamental research on the topic
``Mathematical
Foundations of the Theory of Error-Correcting Codes" and
was also
partly supported by the Higher School of Modern
Mathematics MIPT (project No. FSMG-2024-0048), Moscow College of Physics
and Technology, 1 Klimentovskiy per., Moscow, Russia.

%We thank to the anonymous referees for their remarks and corrections, which have enabled us to greatly improve the paper.
\bibliographystyle{elsarticle-harv}
\bibliography{crmds.bib}

\end{document}